\documentclass[12pt]{article}

\usepackage{amssymb}
\usepackage{amsbsy}
\usepackage{graphicx}
\usepackage{amsmath}
\usepackage{times}

\begin{document}

\title{Razborov flag algebras as algebras of measurable functions}
\author{Tim Austin}
\date{}

\maketitle

% ENVIRONMENTS & THEOREMS

\newenvironment{nmath}{\begin{center}\begin{math}}{\end{math}\end{center}}

\newtheorem{thm}{Theorem}[section]
\newtheorem{lem}[thm]{Lemma}
\newtheorem{prop}[thm]{Proposition}
\newtheorem{cor}[thm]{Corollary}
\newtheorem{conj}[thm]{Conjecture}
\newtheorem{dfn}[thm]{Definition}
\newtheorem{prob}[thm]{Problem}

% COMMANDS

\newcommand{\A}{\mathcal{A}}
\newcommand{\B}{\mathcal{B}}
\newcommand{\K}{\mathcal{K}}
\renewcommand{\Pr}{\mathrm{Pr}}
\newcommand{\s}{\sigma}
\renewcommand{\P}{\mathcal{P}}
\newcommand{\F}{\mathcal{F}}
\newcommand{\M}{\mathcal{M}}
\renewcommand{\O}{\Omega}
\renewcommand{\o}{\omega}
\renewcommand{\S}{\Sigma}
\newcommand{\T}{\mathrm{T}}
\newcommand{\co}{\mathrm{co}}
\newcommand{\e}{\mathrm{e}}
\renewcommand{\d}{\mathrm{d}}
\renewcommand{\i}{\mathrm{i}}
\renewcommand{\l}{\lambda}
\newcommand{\U}{\mathcal{U}}
\newcommand{\G}{\Gamma}
\newcommand{\g}{\gamma}
\renewcommand{\L}{\Lambda}
\newcommand{\hcf}{\mathrm{hcf}}
\renewcommand{\a}{\alpha}
\newcommand{\bbN}{\mathbb{N}}
\newcommand{\bbR}{\mathbb{R}}
\newcommand{\bbT}{\mathbb{T}}
\newcommand{\bbZ}{\mathbb{Z}}
\newcommand{\bbE}{\mathbb{E}}
\newcommand{\Sym}{\mathrm{Sym}}

\newcommand{\bb}[1]{\mathbb{#1}}
\renewcommand{\rm}[1]{\mathrm{#1}}
\renewcommand{\cal}[1]{\mathcal{#1}}

\newcommand{\para}[1]{\paragraph{#1}}
\newcommand{\qed}{\nolinebreak\hspace{\stretch{1}}$\Box$}
\newcommand{\fin}{\nolinebreak\hspace{\stretch{1}}$\lhd$}

\begin{abstract}
These are some brief notes on the translation from Razborov's
recently-developed notion of flag algebra (\cite{Raz07}) into the
lexicon of functions and measures on certain abstract Cantor spaces
(totally disconnected compact metric spaces).
\end{abstract}

\section{The objects of interest}

Consider a universal first-order theory $T$ with equality in a
language $L$ that contains only predicate symbols; assume $T$ has
infinite models.  Examples include the theories of undirected and
directed graphs and hypergraphs, possibly with loops.

In~\cite{Raz07}, Razborov develops a formalism for handling the
`leading order' statistics of large finite models of such theories.
The central objects of his theory are the positive
$\bbR$-homomorphisms of `flag algebras'. Here we shall relate these
to measures on a subset of the Cantor space $\prod_{i \leq
k}K_i^{\bbN^i}$ (where each $K_i$ is itself some Cantor space) which
is compact and such that both set and measure are invariant under
the canonical coordinate-permuting action of $\Sym_0(\bbN)$.
Examples of such Cantor spaces include the spaces of models over
$\bbN$ of certain kinds of theory, and our identification will begin
here. In particular we will identify (under some restrictions) a
flag algebra with an algebra of measurable functions on the
underlying Cantor space of models, and this will lead to an
identification of the positive homomorphisms of the flag algebra
with certain measures by virtue of classical results of functional
analysis.

This will take us through the first three sections of~\cite{Raz07}.
Section 4, 5 and 6 of~\cite{Raz07} relate to a more precise
variational analysis of certain examples of these homomorphisms, and
we will not discuss this here. For certain special examples of the
theories appearing in the study of flag algebras (particularly the
theories of hypergraphs and directed hypergraphs), the associated
probability measures on $\prod_{i \leq k}K_i^{\bbN^i}$ fall into the
classical study of `partial exchangeability' in probability theory,
and quite complete structure theorems describing all such measures
are available following work of Hoover~\cite{Hoo79,Hoo82},
Aldous~\cite{Ald81,Ald82,Ald85} and Kallenberg~\cite{Kal92}.
However, these general results do not encompass the variational
analysis undertaken in Section 5 of~\cite{Raz07} and
in~\cite{Raz07.2}, and it is not clear that adopting the older
probabilistic formalism makes much difference to this.  We direct
the reader to~\cite{Aus--ERH} for a survey of this older theory and
its relations with combinatorics.

Although it is not assumed in~\cite{Raz07}, let us here assume for
simplicity that the predicates of $T$ have bounded maximum arity,
say by $k$. The general case can be recovered as a suitable `inverse
limit' in whichever picture is chosen for the description that
follows; the necessary modifications are routine. We will also
assume that $L$ is countable, although this assumption can also be
dropped with only a cosmetic increase in complexity, and will write
$L_i$ for the set of predicates in $L$ of arity $i$ for $i \leq k$.
Otherwise we adopt the notation and definitions of \cite{Raz07}, and
have tried to avoid conflict with such further notation as we
introduce ourselves.

Given this arity bound, a general countably infinite model of such a
theory $T$ can be identified with a non-uniform coloured looped
directed hypergraph of maximal rank $k$, in which for each rank $i$
the `colour' of an $i$-tuple $a \in \bbN^{i}$ is a truth-assignment
for $\ell(a)$ for each predicate $\ell \in L$ of arity $i$, and so
may be interpreted as a point of $K_i := \{0,1\}^{L_i}$.   Let us
work henceforth only with models whose set of variables (or
`vertices') is a subset of $\bbN$. Thus we shall view out models of
$T$ as points of the Cantor space $\prod_{i \leq k}K_i^{\bbN^i}$
satisfying all additional constraints imposed by the sentences of
$T$; we shall denote the subset of these points by $\O$. In general,
given a model $M$ of $T$ with vertex set $S$ and $S_1\subseteq S$ we
shall write $M|_{S_1}$ for the submodel of $M$ with vertex set
$S_1$. If $T$ is \emph{free} then all points of this space are
possible; otherwise the constraints imposed by $T$ carve out some
subspace of $\prod_{i \leq k}K_i^{\bbN^i}$, which is an intersection
of a family of finite-dimensional subsets corresponding to
individual constraints (since any individual interpretation of a
sentence in $T$ over some particular finite set of vertices in
$\bbN$ simply carves out some clopen subset of $\prod_{i \leq
k}K_i^{\bbN^i}$ depending only on those vertices as coordinates).
This set of models is therefore a compact subspace of an abstract
Cantor space, and so itself an abstract Cantor space. This subspace
is also clearly invariant under permutations of the variable set
$\bbN$; let us write $\Sym_0$ for the group of finitely-supported
such permutations (it will prove convenient to have imposed the
finite-support condition when we come to work with this group
later). For convenience, if a theory $T$ is clear from the context,
we shall refer to a \textbf{cylinder subset} of $\O$ to mean the
intersection of $\O$ with a cylinder subset of the big product space
$\prod_{i \leq k}K_i^{\bbN^i}$.

In fact, our analytic formalism will apply to any closed
$\Sym_0$-invariant subspace of $\prod_{i\leq k}K_i^{\bbN^i}$ for
compact metric spaces $K_1$, $K_2$, \ldots, $K_k$, and the reader is
free to view $\O$ in this generality.

The relevant data for what follows are the nonempty
$\Sym_0$-invariant compact subset $\O$ together with the canonical
action of the finitely-supported permutation group $\Sym_0$ on the
vertex set $\bbN$, and with its usual topological structure and
Borel $\s$-algebra $\S$.  By mild notational abuse we will identify
an element $g \in \Sym_0$ directly with the corresponding
autohomeomorphism of $\O$.

We need to allow a few more constructions. Henceforth all finite
models of $T$ will be assumed to have vertex sets equal to initial
segments of the integers.  For a given such finite model, say $\s$
on the vertex set $[k]$, we let $\O_\s$ be the subset of those $\o
\in \O$ with $\o|_{[k]} = \s$, and let $\S_\s$ be its Borel
$\s$-algebra (which is just the ideal in $\S$ of the subsets of
$\O_\s \subseteq \O$ in the usual way). This $\s$ is called a `type'
in~\cite{Raz07}. We write $0$ for the empty type. We now consider
isomorphism classes of those finite $T$-models that contain a
distinguished copy of $\s$ as a submodel; that is, of pairs $F =
(M,\theta)$ for $M$ a finite model and $\theta$ a particular
embedding $\s \hookrightarrow M$. These extensions are called
`$\s$-flags' in~\cite{Raz07}, and $\F^\s$ is written for the
collection of them, $\F^\s_\ell$ for the subcollection of those with
vertex set of size $\ell$ (so $\F^\s_\ell = \emptyset$ unless $\ell
\geq k$), and we can now specify the obvious notions of embedding
and isomorphism for flags.

Restricting our attention to $\s$-flags $F = (M,\s)$ with $M$ a
model on some $[\ell]$ and $\s$ itself serving as its distinguished
copy in $M$ (so we have implicitly ordered the vertices of $M$ so
that $\s = M|_{[k]}$), we can identify $F$ with a finite-dimensional
cylinder set $A_F$ in $\O$ contained in $\O_\s$:
\[A_F := \{\o \in \O:\ \o|_{[\ell]} = M\}\]
(of course, this is also just $\O_M$ in our earlier notation
associating $\O_\s$ to $\s$; our choice of different letters
reflects the different r\^oles of $M$ and $\s$ here).  In fact, we
could have made any assignment $F \mapsto A_F$ as above but
corresponding to any other choice of locations in $\bbN$ for the
vertices of $M$ that are not vertices of $\s$, and it would not
matter; the above choice is conveniently concrete.

Occasionally we shall need to relate flags over different types $\s$
and the associated spaces $\O_\s$.  In general, the appearance of
$\s$ as a submodel of $\s'$ does not guarantee that specifically $\s
= \s'|_{[k]}$ for some $k$; but we can always identify $\O_{\s'}$ as
a subspaces of $\O_\s$ by suitably reordering the vertices of $\s'$.
In~\cite{Raz07} this technical matter is more-or-less avoided
altogether owing to the early decision to work entirely with
isomorphism classes of models; in the picture of the spaces $\O_\s$,
however, we pay this modest price for the sake of retaining a more
concrete picture of the points of this space and (eventually)
measures on them, which have more classical and easily-analyzed
structures in certain other respects.

Finally, we will want to consider coordinate-permutations that
preserve $\O_\s$; that is, that permute only coordinates in
$\bbN\setminus[k]$.  Let us call the group of these $\Sym_\s \leq
\Sym_0$ (even though it actually depends only on $k$).

\section{Some analysis and measure theory}

Having laid out the objects of study in the previous section, we
will now recall those additional functional-analytic ideas in terms
of which we'll later give an account of flag algebras. The
functional analysis needed does not extend beyond a graduate-level
introduction to probability theory and the contents of any good
second course on functional analysis; Yosida~\cite{Yos95}, for
example, covers our needs.

Let us write $C(\O_\s)$ and $\M(\O_\s)$ for the usual Banach spaces
of real-valued continuous functions and signed Radon measures on
$\O_\s$ respectively; the Riesz-Kakutani representation identifies
$\M(\O_\s)$ isometrically with $C(\O_\s)^\ast$. Moreover, we write
$\M^\s$ for the subspace of measures supported on $\O_\s$ and
invariant under finitely-supported permutations of the coordinates
in $\bbN\setminus[k]$; we will call these
\textbf{$\s$-exchangeable}, and sometimes abbreviate this to just
`exchangeable' when $\s$ is clear from the context. We write
$(\M^\s)^\perp$ for the annihilator of this considered as a subspace
of $C(\O_\s)^\ast$,
\[(\M^\s)^\perp = \{f \in C(\O_\s):\ \langle f,\mu\rangle
= 0\ \forall \mu \in \M^\s\};\] as usual, the dual-of-the-quotient
Banach space $\big(C(\O_\s) / (\M^\s)^\perp\big)^\ast$ can be
isometrically identified with $\M^\s$.  Let $q_\s$ be the quotient
map $C(\O_\s) \to C(\O_\s) / (\M^\s)^\perp$.

Now, given any $f \in C(\O_\s)$, it `looks the same' as any $f\circ
g$ for $g \in \Sym_\s$ to all measures in $\M^\s$. Write $\T_\s$ for
the tail $\s$-subalgebra $\bigcap_{m \geq k+1}\S_{[k]\cup
\{m,m+1,\ldots\}}$ of $\S_\s$, and (with a slight abuse of notation)
$L^\infty(\T_\s)$ for the space of bounded $\T_\s$-measurable
functions that are defined $\mu$-a.e. for every $\mu \in \M^\s$ and
under the equivalence relation of ``equality $\mu$-a.e. for every
$\mu \in \M^\s$''. Clearly these are invariant under the action of
$\Sym_\s$. By the pointwise ergodic theorem for the amenable group
$\Sym_\s$ (or any more elementary argument for this very specialized
example of a $\Sym_\s$-system), we may take the average of the
compositions $f\circ g$ over different $g$ to obtain some
$\T_\s$-measurable function $\bar{f}$ on $\O_\s$ which is defined
$\mu$-almost everywhere for every $\mu \in \M^\s$ and is invariant
under $\Sym_\s$, and hence actually specifies a member of
$L^\infty(\T_\s)$. Observe that $\bar{f} = \bar{h}$ for $f,h \in
C(\O_\s)$ if and only if $f - h \in (\M^\s)^\perp$, and so our map
$f \mapsto \bar{f}$ embeds $C(\O_\s)/(\M^\s)^\perp$ as a subspace
$V^\s$ of $L^\infty(\T_\s)$; general nonsense now shows also that
this is an isometric embedding, so $V^\s$ is closed.

Furthermore, $V^\s$ is actually a sub\emph{algebra} of
$L^\infty(\T_\s)$.  To show that it is closed under multiplication,
we suppose $f,h \in C(\O_\s)$, and now consider the products $f\cdot
(h\circ g)$ for any sequence of permutations $g$ that pushes $h$
`further and further out', in the following sense: for any $m \geq
1$, there are finite subsets $A,B \subset \bbN\setminus [k]$ such
that $f$ and $h$ are uniformly $(1/m)$-close to functions depending
only on vertices in $A$ and $B$ (respectively), and now we choose
$g$ that moves all points of $B$ into $\bbN \setminus A$.  Letting
$m \to \infty$ this gives a sequence $g_m$ for which, in terms of
their dependence on coordinates, $f$ and $h\circ g_m$ are closer and
closer to independent.

Now the point is that the quotients $q_\s(f\cdot h\circ g)$ converge
in $C(\O_\s) / (M^\s)^\perp$ to a member that depends only on
$q_\s(f)$ and $q_\s(h)$ --- this follows from an elementary
step-function approximation argument and use of all the permutation
invariance. Let us call this the \textbf{asymptotic product} of
$q_\s(f)$ and $q_\s(g)$.  It is now a routine exercise to check this
is actually a $\rm{C}^\ast$-algebra product on $C(\O_\s) /
(M^\s)^\perp$ and corresponds exactly to the usual product of
functions in $V^\s$.

(Alternatively, one can find an actual continuous function on
$\O_\s$ whose image in $V^\s$ represents this product; let us
illustrate one cheap way to do this in case $\s = 0$. Let
$\psi_1:2\bbN \to \bbN$ and $\psi_2:2\bbN + 1 \to \bbN$ be
bijections, and let us use the same letters for the corresponding
adjoint maps $\O \to \O^{(2\bbN)}\stackrel{\phi_1}{\cong} \O$ and
$\O \to \O^{(2\bbN+1)}\stackrel{\phi_2}{\cong} \O$, where we
temporarily write $\O^{(S)}$ for the space of models of $T$ with
vertex set $S$. Then clearly the functions $f\circ\phi_1\circ\psi_1$
and $g\circ\phi_2\circ\psi_2$ depend on disjoint sets of
coordinates; their product is the sought representative for
$q_\s(f)q_\s(g)$.)

Thus we have identified $C(\O_\s)/(\M^\s)^\perp$ with a closed
subalgebra $V^\s$ of $L^\infty(\T_\s)$ (with a newly-defined
product). Let us call functions in $V^\s$ \textbf{simple} if they
are the images of simple (equivalently, finite-dimensional)
functions in $C(\O_\s)$. We can describe the simple functions
naturally as follows: to any fixed nonempty cylinder set $A
\subseteq \O_\s$ depending on coordinates in $J \subset
\bbN\setminus [k]$ corresponds a collection of finite models of the
theory $T$ on the vertex set $J$ (with some multiplicities), and now
the averaged-function $\overline{1_A}(\o)$ for $\o \in \O_\s$ is
just the sum of the densities with which each of those finite models
appears isomorphically as a submodel of $\o$ (now summing over the
multiplicities). Referring to such a function $\overline{1_A}$ for
$A$ corresponding to a single model on $J$ (so that our general $A$
is a disjoint union of such) as a \textbf{statistics function}, the
simple functions in $V^\s$ are now just linear combinations of
statistics functions. We write $V^\s_0$ for the dense subspace of
these.

\section{Description of flag algebras}

We will now identify a flag algebra and its homomorphisms with a
family of measurable functions and an associated set of probability
measures, and show how various results pertaining to the former
translate into facts about the latter.  In general we will refer to
the two resulting pictures as the `flag algebra picture' and the
`measure theoretic picture' respectively.  We will follow the
sectional structure of Sections 2 and 3 of~\cite{Raz07}. We will not
broach the more specific optimization problem studied in Sections 4
and 5 (or, for that matter, the more complete result
of~\cite{Raz07.2}); while they too can presumably be translated into
the measure theoretic picture, aside from rendering the underlying
objects in a better-established analytic light this doesn't seem
greatly to change the arguments that are involved.

Informally, the background purpose of~\cite{Raz07} is to set up a
family of `proxies' for the statistics of large models of $T$ that
enjoy some additional `analytic' structure making them easier to
handle for the study of certain kinds of question; two such
questions on extremal statistics are then analyzed in these terms in
Section 5 of~\cite{Raz07} and in the follow-up paper~\cite{Raz07.2}.
In~\cite{Raz07} these proxies are certain $\bbR$-valued
homomorphisms of flag algebras, and the actual elements of the flag
algebras are of secondary importance (although they do continue to
appear in a supporting role occasionally later in the paper).  The
manipulation of flag algebras is mostly for the purpose of setting
up these homomorphisms. However, in the measure theoretic picture we
can say at once what these homomorphisms correspond to --- they are
the ergodic $\Sym_\s$-invariant probability measures on $\O_\s$ ---
and so the effort we expend on setting up the correspondence between
flag algebras and certain spaces of measurable functions will
ultimately be required only to show that the homomorphisms of the
former really are identified with a priori-known objects related to
the latter.

\subsubsection*{Section 2}

In~\cite{Raz07}, $\bbR\F^\s$ denotes the free $\bbR$-vectorspace on
$\F^\s$. We have already selected for each $F \in \F^\s$ some
associated clopen subset $A_F\subseteq \O_\s$, and we will now
extend this by associating (with some careful normalization) to each
member $\sum_j\l_kF_j \in \bbR\F^\s$ a linear combination of the
indicator functions $1_{A_F}$ of these $A_F$, each of which is a
continuous function since $A_F$ is clopen. This will define a linear
operator $\Phi:\bbR\F^\s \to C(\O_\s)$ with image some funny-shaped
subspace contained within the space of simple functions in
$C(\O_\s)$. It turns out that in order to make contact between the
calculus of~\cite{Raz07} and addition and multiplication of
functions on $\O_\s$ we need to introduce some nontrivial
normalizing constants: our final identification is
\begin{eqnarray*}
F &\mapsto& \frac{1}{|F|!}1_{A_F};\\
\sum_j\l_kF_j &\mapsto& \sum_j\frac{\l_k}{|F_j|!}1_{A_{F_j}}.
\end{eqnarray*}

The arbitrariness in our choice of subset $A_F$ corresponding to $F$
is reflected in a similar arbitrariness in this linear map $\Phi$
into $C(\O_\s)$; however, this disappears at the next step, when we
define a flag algebra as a quotient of $\bbR\F^\s$.

Let $\K^\s$ denote the subspace of $\bbR\F^\s$ generated by the
linear combinations $\tilde{F} - \sum_{F \in
\F^\s_\ell}p(\tilde{F},F)F$ for different $\ell \geq
|V(\tilde{F})|$, and for certain real numbers $p(\tilde{F},F)$ given
in Definition 1 of~\cite{Raz07} (we will suppress their exact form
here). From that definition one can check at once that, given our
chosen normalization in the definition of $\Phi$ above, the values
$p(\tilde{F},F)$ are such that $\Phi(\K^\s)$ is precisely the
following subspace of $\Phi(\bbR\F^\s)$:
\begin{quote}
$\Phi(a) \in \Phi(\K^\s)$ if and only if, as a simple function on
$\O_\s$, $\Phi(a)$ may be chopped up further into a linear
combination of simple functions corresponding to cylinder sets over
some common large finite subset of $\bbN$, say $\Phi(a) = \sum_{t\in
J}b_t1_{B_t}$, such that we can cluster this sum according to some
partition $J = \bigcup_{i\in I}J_i$,
\[\Phi(a) = \sum_tb_t1_{B_t} = \sum_{i \in I}\sum_{t \in J_i}b_t1_{B_t},\]
so that for each $i \in I$:
\begin{itemize}
\item all the $B_t$ for $t \in J_i$ are isomorphic to some fixed
$A_{F_i}$, $F_i \in \F^\s$;
\item and $\sum_{t\in J_i}b_t = 0$.
\end{itemize}
\end{quote}
Alternatively, we can describe this by saying that some
re-arrangement of the terms $b_t1_{B_t}$ by different permutations
of $\bbN\setminus[k]$ sums to zero: that is, there are $g_1$,
\ldots, $g_t \in \Sym_\s$ with $\sum_tb_t1_{g_t(B_t)}$ exactly
canceling to zero.

It now follows from a little compactness argument that $\K^\s$ is
also precisely the set of those $a \in \bbR\F^\s$ for which $\Phi(a)
\in C(\O_\s)$ is in the annihilator $(\M^\s)^\perp$ (indeed, that
$\Phi(\K^\s) \subseteq (\M^\s)^\perp$ is immediate; for the opposite
inclusion we can argue that if $a \not\in \K^\s$ then by witnessing
our inability to find a decomposition of $\Phi(a)$ as above for
cylinder sets over larger and larger finite subset of $\bbN$, we can
extract some sequence of members of $\M(\O_\s)$ that are invariant
for the vertex-permutations in some corresponding exhausting
sequence of finite subsets of $\Sym_\s$\footnote{Notice that here is
an appeal to our restriction to finitely-supported permutations.},
and that converge to some member of $\M^\s$ that does not annihilate
$\Phi(a)$).

We now consider the quotient space $\A^\s := \bbR\F^\s/\K^\s$. By
the above, $q_\s\circ\Phi$ factors through this quotient to give an
injective map $\Psi:\A^\s \to V^\s$.  Moreover, since any finite
dimensional cylinder set contained in $\O_\s$ is equivalent under
$\M^\s$ to some $A_F$ upon a suitable permutation of coordinates in
$\bbN\setminus[k]$, the image of $\Psi$ is actually the subspace
$V^\s_0$ of \emph{all} simple functions in $V^\s$; as such, it is
dense.

The next step is to consider products.  The definition of a product
for $a,b \in \bbR\F^\s/\A^\s$ given in~\cite{Raz07} now just
translates into the product of $\Psi(a)$ and $\Psi(b)$ as
$L^\infty$-functions on $\O_\s$; and the proof in~\cite{Raz07} that
this product is well-defined mostly becomes a proof in the measure
theoretic picture that this product remains in the image of $\Psi$
(that is, in $V^\s_0$). This completes our identification of the
flag algebra $\A^\s$ with the dense subalgebra $V^\s_0$ of $V^\s$,
which is itself a norm-closed Banach subalgebra of
$L^\infty(\T_\s)$.  The basic properties of the product contained in
Lemma 2.4, for example, are now immediate.

In Subsection 2.2 of~\cite{Raz07} is introduced the `downward
operator'. This applies when we have a submodel $\s'$ of $\s$, say
with vertex set $[k']$ for some $k' \leq k$. By re-labeling the
vertices of $\s'$ if necessary, we can identify $\O_\s$ with a
clopen subset of $\O_{\s'}$ (with $A_{(\s',\s)}$, in fact); since
this subset is clopen, the extension operator $J:C(\O_\s)
\hookrightarrow C(\O_{\s'})$ obtained by extending a continuous
function on $\O_\s$ to be identically zero elsewhere is well-defined
(in particular, its output is still a continuous function).  One
checks at once that $J$ factorizes through the quotients $C(\O_{\s})
\to V^\s \subseteq L^\infty(\T_{\s})$ and $C(\O_{\s'}) \to V^{\s'}
\subseteq L^\infty(\T_{\s'})$.  This factorization is the measure
theoretic picture of the downward operator. Once again, a
normalizing factor appears in~\cite{Raz07} to make the sums come out
right.

Subsection 2.3 turns to the `upwards operator'. The definition of
this and the properties it enjoys depend more heavily on the precise
shape of the theory $T$ than most of the foregoing, and we shall not
translate the results of this subsection in detail. Instead, let us
examine only a leading special case of this operator.

Given again some type $\s$ on $[k]$ extending $\s'$ on $[k']$, $k'
\leq k$, we now wish to make a passage from $V^{\s'}$ to $V^\s$ in
the following way. Any point of $\O_\s$ defines a point of
$\O_{\s'}$ simply by ignoring the vertices in
$\{k'+1,k'+2,\ldots,k\}$ (and so sending $\bbN\setminus [k]$ to
$(\bbN-(k-k'))\setminus [k] = \bbN\setminus [k']$); this gives a
continuous map from $\O_{\s}$ to $\O_{\s'}$, so that composition
gives a homomorphism $C(\O_{\s'}) \to C(\O_{\s})$. This \emph{may}
descend to a map $V^{\s'} \to V^\s$ which is then necessarily also a
homomorphism: this requires that any member of $(\M^{\s'})^\perp$ be
sent to a member of $(\M^\s)^\perp$ by this composition map, which
in some sense tells us that the space of models $\O_{\s}$ is still
`large enough' to support a sufficiently large collection of
$\Sym_{\s}$-invariant measures compared with $\O_{\s'}$.  A formal
version of this property appears in the flag algebra picture as a
condition that a certain member of $\A^\s$ is not a zero-divisor,
and under this assumption the existence of a suitable homomorphism
(which translates into the abovementioned factorizability) is proved
directly for flag algebras as Theorem 2.6 in~\cite{Raz07}; this
property is also related to the more immediate property of a theory
that it have the `amalgamation property'.

This is an example that can be extended to relate the algebras
$\A^\s$ (or $V^\s$) arising from two different theories $T_1$, $T_2$
given an interpretation of one in terms of the other: this
interpretation again defines a continuous map from the Cantor space
corresponding to one to that corresponding to the other, and again
we then face the question of whether the resulting homomorphism of
algebras of continuous functions descends under our quotienting
operation. The considerations to this purpose in~\cite{Raz07} apply
in this general setting, but we will not examine the details further
here.

When a homomorphism can be obtained, some of the other results that
follow on the properties of this map are now translations of certain
basic facts for concrete spaces of measurable functions (to which
they already appear structurally very similar): Theorem 2.8(a)
of~\cite{Raz07}, for example, asserts in our picture that
$\bbE_\mu[fh\,|\,\Xi] = f\bbE_\mu[h\,|\,\Xi]$ if $f$ is already
$\Xi$-measurable, and 2.8(b) is the rule of iterated conditional
expectations.

\subsubsection*{Section 3}

The overall approach of the first three sections of~\cite{Raz07} is
to define a flag algebra \emph{first}, and then to obtain a
collection of `proxies' or `limit objects' for the statistics of
large models of a theory in terms of them.  The r\^ole of these is
to be played by the multiplicative functions $\phi:\A^\s \to \bbR$
that are non-negative on the image of any single flag $F \in \F^\s$;
the set of these is written $\rm{Hom}^+(\A^\s,\bbR)$.

Having identified the flag algebra $\A^\s$ as the dense subalgebra
$V^\s_0$ of the commutative C$^\ast$-algebra $V^\s$ so that the
images of single flags correspond to the single statistics
functions, we can easily check that non-negativity on these implies
non-negativity on any member of $V^\s_0$ that is itself a
non-negative-valued function; this follows easily since a
non-negative simple function can always be written as a linear
combination of indicator functions with non-negative coefficients.

We now observe that \emph{given the non-negativity} of such a $\phi$
it can be extended to a multiplicative linear functional on the
whole of $V^\s$. We can now, if we wish, apply certain standard
representation theorems to this space: either by further exploiting
its vector lattice structure following the results of Yosida and
Kakutani, as presented in Section XII.5 of~\cite{Yos95}, or by
complexifying the construction so far and using the (arguably more
popular) Gelfand-Naimark Theorem. At any rate, this identifies
$\phi$ with a point of the spectrum of $V^\s$ (for one or other
interpretation of `spectrum'). In fact, this identification is
more-or-less implicit in Remark 4 of Subsection 3.2 of~\cite{Raz07},
although there we still require some of the basic structure of
$\rm{Hom}^+(\A^\s,\bbR)$ to have been idenfitied.

However, given our identification with the measure-theoretic
picture, we have an alternative to the above. Since $V^\s \cong
C(\O_\s)/(\M^\s)^\perp$, as a linear functional on $V^\s$ we can
identify $\phi$ with a member of $\M^\s$ (uniquely, since
$(C(\O_\s)/(\M^\s)^\perp)^\ast \cong \M^\s$): a $\Sym_\s$-invariant
measure on the space $\O_\s$ that we started with, rather than a
point of some abstractly-produced new space $\rm{Spec}\,V^\s$.  This
has the advantage that we retain the theory itself $T$, coded as it
is into the `shape' of the space $\O_\s$. It is now easy to check
that those measures in $\M^\s$ that are \emph{multiplicative} on
$V^\s$ are precisely the ergodic $\Sym_\s$-invariant probability
measures on $\O_\s$ (these are multiplicative on $V^\s$ since any
member of $V^\s$ is $\mu$-a.s. constant if $\mu \in \M^\s$ is
ergodic). This establishes the identification of the flag algebra
homomorphisms with exchangeable measures~\footnote{Although we note
in passing that obtaining them from a prior construction of the
function algebras is reminiscent of the alternative route into the
theory of integration and measure established by Segal
in~\cite{Seg54}; an approach that has had lasting consequences for
the formulation of `non-commutative integration' in the setting of
general von Neumann algebras.}.

As remarked in Remark 3 of Section 3 of~\cite{Raz07}, working with
arbitrary multiplicative linear functionals on $\A^\s$ is
problematic: the point is that without non-negativity these need not
be extendable to the whole of $V^\s$ at all (equivalently, they may
not be continuous for the norm topology of $V^\s_0$).  Indeed, the
fact, mentioned in Remark 3, that $\A^\s$ is a
non-finitely-generated free commutative algebra over $\bbR$, is
precisely what would allow us to construct a discontinuous such
functional using the axiom of choice.

Now the order defined on $\A^\s$ in Definition 5 of Section 3
of~\cite{Raz07} (in terms of the above notion of `positivity' for a
homomorphism, which must be introduced first) is precisely the usual
pointwise order on $V^\s_0$ as a set of real-valued functions; in
the setting of abstract flag algebras, where we are unable to define
anything `pointwise', the functionals of $\rm{Hom}^+(\A^\s,\bbR)$
are needed as a replacement to formulate this definition.  Given
this, Theorem 3.1 requires only that composition with a
homeomorphism and conditional expectation are non-negative operators
between function spaces.  Also, the `probability-like' convergence
results of Subsection 3.1 (somewhat based on~\cite{LovSze06}) now
really are about the classical vague topology on a set of
probability measures.

In Subsection 3.2 of~\cite{Raz07} averages of homomorphisms are
taken with respect to actual probability measures; in the
measure-theoretic picture these become classical ergodic
decompositions. Specifically, given some extension $\s$ of $\s_0$,
we know that upon ordering the vertices of $\s$ so that $\s_0 =
\s|_{[k_0]}$, then $\O_\s = A_{(\s_0,\s)}$ becomes a clopen subset
of $\O_{\s_0}$; and now if $\mu$ is an ergodic
$\Sym_{\s_0}$-invariant probability measure on $\O_{\s_0}$ for which
$\mu(\O_\s) > 0$ we may consider the conditioned measure
$\mu(\,\cdot\,|\O_\s) := \mu(\,\cdot\,\cap\O_\s)/\mu(\O_\s)$
(defined free from any measure-theoretic ambiguity, since
$\mu(\O_\s) > 0$).  Since $\Sym_\s \leq \Sym_{\s_0}$ fixes $\O_\s$,
this defines now a $\Sym_{\s}$-invariant probability measure on
$\O_\s$; however, it may not be ergodic under the action of this
subgroup $\Sym_{\s}$, and the resulting `ensemble' of homomorphisms
obtained in this subsection is simply its ergodic decomposition.

After translation, most of the results of this subsection are
identified with the standard results for actual probability measures
that they mimic. For example, the argument that a suitable
probability measure on $\rm{Hom}^+(\A^\s,\bbR)$ exists with
barycentre a given member of $\rm{Hom}^+(\A^{\s'},\bbR)$ for
Definition 8 now asserts the existence of ergodic decompositions.

Finally, the results of Subsection 3.3 become ordinary inequalities
and continuity results for functions and measures. Theorem 3.14 is
the conditional Cauchy-Schwartz inequality.  Theorems 3.17 and 3.18
relate to iterated conditional expectations and the pointwise order
of functions; essentially the same proofs as in~\cite{Raz07} are now
the proofs of the basic facts about measures and functions.

Theorems 3.15 and 3.16 are more specific to the study of models of a
theory $T$, and here the matter of which formalism we choose for
their proof seems quite unimportant; we forego giving the
measure-theoretic details.

\bibliographystyle{abbrv}
\bibliography{Razborovm}

\vspace{10pt}

\parskip 0pt

\textsc{Department of Mathematics\\ University of California at Los
Angeles, Los Angeles, CA 90095-1555, USA}

Email: \verb|timaustin@math.ucla.edu|

Web: \verb|http://www.math.ucla.edu/~timaustin|

\vspace{7pt}

%DATE!

\end{document}